\documentclass{amsart}

\usepackage{amssymb}
\usepackage[english]{babel}
\usepackage[dvips]{graphicx}

\usepackage[numeric, short-months]{amsrefs}

\newtheorem{theorem}{Theorem}[section]
\newtheorem{corollary}[theorem]{Corollary}
\newtheorem{proposition}[theorem]{Proposition}

\theoremstyle{definition}
\newtheorem{definition}[theorem]{Definition}

\newcommand{\Alphabet}[1]{ {X} _{#1} }
\newcommand{\InfiniteSet}[1][m]{ { \Alphabet{#1} ^ \omega } }
\newcommand{\FiniteSet}[1][m]{ { \Alphabet{#1} ^ \ast } }
\newcommand{\States}[1]{ {Q} _{#1} }

\newcommand{\AutomataSet}[2]{ \Automaton[{#1 \times #2}] }
\newcommand{\Automaton}[1][]{ { {A} _{#1} } }
\newcommand{\AExample}[1]{ {B} _{#1} }
\newcommand{\MAutomaton}[3]{ \left( \Alphabet{#2}, \States{#1}, \pi _{#3}, \lambda _{#3} \right) }

\newcommand{\Semigroup}[1]{ {S} _{#1} }
\newcommand{\ATMonoid}[1]{ \left\langle {\begin{array}{*{20}c}
   {e, f_0, f_1} & \vline &  {#1} \\
\end{array}} \right\rangle }
\newcommand{\ATSemigroup}[1]{ \left\langle {\begin{array}{*{20}c}
   {f_0, f_1} & \vline & {#1} \\
\end{array}} \right\rangle }
\newcommand{\ATSemigroupThree}[1]{ \left\langle {\begin{array}{*{20}c}
   {f_0, f_1, f_2} & \vline &  {#1} \\
\end{array}} \right\rangle }

\newcommand{\SequenceIndex}[2]{ { \left\{ {#1 _{m}}, \, m \geq {#2} \right\} } }

\newcommand{\Growth}[1]{ \gamma _{#1} }
\newcommand{\GrowthOrder}[1]{ { \left[ {#1} \right] } }
\newcommand{\GrowthAutomaton}[1]{ \Growth{ \Automaton[#1] } }
\newcommand{\GrowthSemigroup}[1]{ \Growth{ \Semigroup{#1} } }
\newcommand{\SphericalGrowthSemigroup}[1]{ \mathord{ \buildrel{ \lower 3pt
\hbox{$\scriptscriptstyle\frown$}} \over {\Growth{}} } _{ \Semigroup{#1} } }
\newcommand{\WordGrowthSemigroup}[1]{ \delta _{ \Semigroup{#1} } }

\newcommand{\n}[1][n]{ ({#1}) }
\newcommand{\Divider}[1] { {\left[ \frac{#1}{2} \right]} }
\newcommand{\Binomial}[2]{\left( {\begin{array}{*{20}c} {#2} \\ {#1} \\ \end{array}} \right)}

\newcommand{\q}{ \mathsf{q} }
\newcommand{\s}{ \mathsf{s} }

\newcommand{\Natural}{ \mathbb{ N } }

\newcommand{\PositiveInteger}{ n \in \Natural }

\begin {document}

\title[On composite and non-monotonic growth of Mealy automata]{On composite
and non-monotonic growth functions of Mealy automata}

\author[I. Reznykov]{Illya I. Reznykov}
\address{Faculty of Mathematics and Mechanics\\Kiev Taras Shevchenko National
University\\vul. Volodymyrska, 64\\Kiev\\Ukraine 01033}
\curraddr{Faculty of Mathematics and Mechanics\\Kiev Taras Shevchenko National
University\\vul. Volodymyrska, 64\\Kiev\\Ukraine 01033}
\email{Illya.Reznykov@ikc5.com.ua}


\date{2004/06/28}


\begin{abstract}
We introduce the notion of composite growth function and provide examples that
illustrate the primary properties of these growth functions. There are provided
examples of Mealy automata that have composite non-monotonic growth functions
of the polynomial growth order. There are described examples of Mealy automata
that have composite monotonic growth functions of intermediate and exponential
growth. The questions concerning the interrelation between the notions
``composite'' and ``non-monotonic'' of a Mealy automaton growth function are
formulated.
\end{abstract}

\maketitle

\section{Introduction}

\par The notion of growth was introduced in the middle of last century
\citelist{\cite{Svarc1955-English} \cite{Milnor1968-Growth}} and was applied to
various geometrical, topological and algebraic objects
\citelist{\cite{Babenko1986-English} \cite{Ufnarovskiy1990-English}}. Mainly,
growth functions of studied objects are non-decreasing monotonic functions of a
natural argument \cite{Babenko1986-English}. For example, the growth function
of a semigroup (group) at the point $n$, $n \ge 0$, equals a number of
different semigroup elements of length $n$. Obviously, the growth function of
an arbitrary semigroup is a non-decreasing monotonic function.

\par Growth of Mealy automata have been studied since the 80th of 20th century
\citelist{\cite{Gecseg1986} \cite{Grigorchuk1988-English}}, and it is close
interrelated with growth of automatic transformation semigroups (groups),
defined by them \cite{Grigorchuk1988-English}. However, the growth functions of
the Mealy automaton and the corresponding semigroup have different properties;
for example, they may have different growth orders
\cite{Reznykov2003-Polynomial}. In the paper we consider special type of growth
functions of Mealy automata --- composite growth functions.

\par Composite function is a function such that it can be described by different
expressions on infinite non-overlapped intervals. There exist Mealy automata
that have composite growth functions of various growth orders. Moreover, some
of these automata have non-monotonic growth functions. There were not known
Mealy automata such that they have non-monotonic growth functions (see survey
\cite{GrigorchukNekrashevichSushchansky2000-English}, and
\cite{ReznykovSushchansky2002-Fibonacci}, \cite{Reznykov2002-Ph.D.-English},
etc.).

\par Preliminaries of the theory of Mealy automata are listed in
Section~\ref{sect:preliminaries}. The notion of composite function is
introduced in Section~\ref{sect:composite_function}. In
Section~\ref{sect:non-monotonic} we provide several examples of Mealy automata,
that have non-monotonic growth functions of the polynomial growth order. In
addition, the Mealy automata with composite growth functions such that one of
its finite differences consists of doubled values, are provided in
Section~\ref{sect:doubled_differnces}. The theorems concerning the main
properties of these automata are formulated, and we list these theorems without
proofs. They can be proved by using the technique, similar to the technique of
\cite{Reznykov2002-Ph.D.-English} (see also \cite{Reznykov2003-Polynomial}). We
are planning to publish the proofs of these theorems in subsequent papers. For
convenience, the propositions, where the normal form of semigroup elements are
formulated, are provided for the most complex of the considered automata.
Moreover, the questions concerning the composite growth functions are appeared,
and some of them are mentioned in Section~\ref{sect:final_remarks}.

\par I would like to thank Igor F. Reznykov and Alexandr N. Movchan, who
helped to find the automata, that are considered in the paper.

\section{Preliminaries} \label{sect:preliminaries}

\par The basic notions of the theory of Mealy automata and the theory of
semigroups can be found in many books, for example
\citelist{\cite{Glushkov1961-English} \cite{Gill1963} \cite{Lallement1979}
\cite{Gecseg1986}}. We use definitions from \cite{Reznykov2003-Polynomial}.

\subsection{Mealy automata}

\par Let's denote the set of all finite words over $\Alphabet{m}$, including
the empty word $\varepsilon$, by the symbol $\FiniteSet$, and denote the set of
all infinite (to right) words by the symbol $\InfiniteSet$. We write  a
function $\phi : \Alphabet{m} \to \Alphabet{m}$ as
\[
    \left( {\begin{array}{*{10}c} {\phi(x_0)} & {\phi(x_1)} &
    \dots & {\phi(x_{m - 1})}\end{array}} \right).
\] Moreover, we have in mind $\Natural = \left\{ 0, 1, 2, \ldots \right\}$.

\par Let $\Automaton = \MAutomaton{n}{m}{}$ be a \textit{non-initial Mealy
automaton} \cite{Mealy1955} with the finite set of states $\States{n} = \left\{
q_0, q_1, \ldots, q_{n - 1} \right\}$, input and output alphabets are the same
and equal $\Alphabet{m}$, $\pi : \Alphabet{m} \times \States{n} \to \States{n}$
and $ \lambda : \Alphabet{m} \times \States{n} \to \Alphabet{m}$ are its
transition and output functions, respectively. The function $\lambda$ can be
extended in a natural way to the mapping $\lambda: \FiniteSet \times \States{n}
\to \FiniteSet$ or to the mapping $\lambda: \InfiniteSet \times \States{n} \to
\InfiniteSet$. The transformation $f _{\q} : \FiniteSet \to \FiniteSet$ ($f
_{\q} : \InfiniteSet \to \InfiniteSet$), defined by the equality \[
    f _{\q} \n[u] = \lambda (u, \q),
\] where $u \in \FiniteSet$ ($u \in \InfiniteSet$), is called
\cite{Glushkov1961-English} \textit{the automatic transformation}, defined by
$\Automaton$ at the state $\q$. The automaton $\Automaton$ defines the set \[F
_{\Automaton} = \left\{ f _{q_0}, f _{q_1 }, \ldots, f _{q_{n - 1}} \right\}\]
of automatic transformations over $\InfiniteSet$. Each automatic transformation
defined by the automaton $\Automaton$ can be written in \textit{the unrolled
form}: \[
    f _{q_i} = \left( {f _{\pi (x_0, q_i)}, f _{\pi (x_1, q_i)},
    \ldots, f _{\pi (x_{m - 1}, q_i)} } \right) \sigma _{q_i},
\] where $i = 0, 1, \ldots, {n - 1}$, and $\sigma _ {q_i}$ is the transformation
over the alphabet $\Alphabet{m}$ defined by the output function $\lambda$: \[
    \sigma _{q_i} = \left( {\begin{array}{*{20}c} {\lambda (x_0, q_i)} &
    {\lambda (x_1, q_i)} & \ldots & {\lambda (x _{m - 1}, q_i)}  \\
\end{array}} \right).
\]

\par Let us define the set of all $n$-state Mealy automata over the $m$-symbol
alphabet by the symbol $\AutomataSet{n}{m}$. The product of Mealy automata is
introduced \cite{Gecseg1986} over the set of automata with the same input and
output alphabet $\Alphabet{m}$ as their sequential applying. Therefore for the
transformations $f _{\q_1, \Automaton[1]}$ and $f _{\q_2, \Automaton[2]}$,
$\q_1 \in \States{n _1}$, $\q_2 \in \States{n _2}$, the unrolled form of the
product $f _{(\q_1, \q_2), {\Automaton[1] \times \Automaton[2]}}$ is defined by
the equality:
\[
    f _{(\q_1, \q_2), {\Automaton[1] \times \Automaton[2]}} = f
    _{\q_1, \Automaton[1]} {f _{\q_2, \Automaton[2]}} = \left( g_0, g_1,
    \ldots, g_{m - 1} \right) \sigma _{\q_1, \Automaton[1]} \sigma _{\q_2,
    \Automaton[2]},
\] where $g _i = f _{\pi_1 ({\sigma _{\q_2, \Automaton[2]}
(x_i), \q_1}), \Automaton[1]} f _{\pi_2 (x_i, \q_2), \Automaton[2]}$, $i = 0,
1, \ldots, {m - 1}$, and all transformations are applied from right to left.

\par The power $\Automaton ^n$ is defined for any automaton $\Automaton$ and
any positive integer $n$. Let us denote $\Automaton ^{\n}$ the minimal Mealy
automaton, equivalent to $\Automaton ^n$. It follows from definition of a
product, that $\left| {\States{\Automaton ^{\n}}} \right| \le \left|
{\States{\Automaton}} \right| ^n$.

\begin{definition} \label{def:growth_automaton} \cite{Grigorchuk1988-English}
The function $\GrowthAutomaton{}$ of a natural argument, defined by \[
    \GrowthAutomaton{} \n = \left| {\States{\Automaton ^{\n}}} \right|, \,
    \PositiveInteger,
\] is called \textit{the growth function}
of the Mealy automaton $\Automaton$.
\end{definition}

\subsection{Semigroups}

\begin{definition} \label{def:transformation_semigroup}
Let $\Automaton = \MAutomaton{n}{m}{}$ be a Mealy automaton. The semigroup
\[
    \Semigroup{\Automaton} = \mathop{sg} \left( f_{q_0}, f_{q_1}, \ldots ,
    f _{q_{n - 1}} \right)
\] is called \textit{the semigroup of automatic transformations,
defined by $\Automaton$}.
\end{definition}

\par Let $\Semigroup{}$ be a semigroup with the finite set of generators $G =
\left\{ s_0, s_1, \ldots, s_{k - 1} \right\}$. The elements of the free
semigroup $G^{+}$ are called \textit{semigroup words} \cite{Lallement1979}. In
the sequel, we identify them with corresponding elements of $\Semigroup{}$.
Let's denote the length of a semigroup element $\s$ by the symbol $\ell (\s)$.

\begin{definition} \label{def:growth_semigroup}
The function $\GrowthSemigroup{}$ of a natural argument such that
\[ \GrowthSemigroup{} \n =
\left| \left\{ {\begin{array}{*{20}c}
  {s \in \Semigroup{}} & \vline & {\ell (s) \le n}  \\
\end{array}} \right\} \right|, \, \PositiveInteger,
\] is called \textit{the growth function of $\Semigroup{}$ relative to the
system $G$ of generators}.
\end{definition}

\begin{definition} \label{def:spherical_growth_semigroup}
The function $\SphericalGrowthSemigroup{}$ of a natural argument such that
\[ \SphericalGrowthSemigroup{} \n =
\left| \left\{ {\begin{array}{*{20}c}
  {s \in \Semigroup{}} & \vline & {s = s_{i_1} s_{i_2} \ldots s_{i_n}, \,
  s_{i_j} \in G, \, 1 \le j \le n}  \\
\end{array}} \right\} \right|, \, \PositiveInteger,
\] is called \textit{the spherical growth function of $\Semigroup{}$ relative
to the system $G$ of generators}.
\end{definition}

\begin{definition} \label{def:word_growth_semigroup}
The function $\WordGrowthSemigroup{}$ of a natural argument such that
\[\WordGrowthSemigroup{} \n =
\left| \left\{ {\begin{array}{*{20}c}
  {s \in \Semigroup{}} & \vline & {\ell (s) = n}  \\
\end{array}} \right\} \right|, \, \PositiveInteger,
\] is called \textit{the word growth function of $\Semigroup{}$
relative to the system $G$ of generators}.
\end{definition}

\par From the definitions~\ref{def:growth_semigroup},
\ref{def:spherical_growth_semigroup} and \ref{def:word_growth_semigroup}, the
following inequalities hold for $\PositiveInteger$:
\begin{equation*} \label{eq:growths_semigroup}
    \WordGrowthSemigroup{} \n \le \SphericalGrowthSemigroup{} \n \le
    \GrowthSemigroup{} \n = \sum \limits _{i = 0} ^n {\WordGrowthSemigroup{}
    \n[i]}.
\end{equation*}
Similarly, from definition~\ref{def:transformation_semigroup} it follows
\cite{Grigorchuk1988-English} that
\[
    \GrowthAutomaton{} \n = \SphericalGrowthSemigroup{\Automaton} \n, \,
    \PositiveInteger.
\]

\subsection{Growth functions}

\par The growth of some object is defined by functions of a natural argument.
One of the most used characteristic of these functions is the notion of growth
order.

\begin{definition} \label{def:growth_order}
Let $\Growth{i} : \Natural \to \Natural$, $i = 1, 2$, are arbitrary functions.
The function $\Growth{1}$ has \textit{no greater growth order} (notation
$\Growth{1} \preceq \Growth{2}$) than the function $\Growth{2}$, if there exist
numbers $C_1, C_2, N_0 \in \Natural$ such that \[\Growth{1} \n \le C_1
\Growth{2} (C_2 n)\] for any $n \ge N_0$.
\end{definition}

\begin{definition} \label{def:growth_equivalence}
The growth functions $\Growth{1}$ and $\Growth{2}$ are equivalent or have
\textit{the same growth order} (notation $\Growth{1} \sim \Growth{2}$), if the
inequalities $\Growth{1} \preceq \Growth{2}$ and $\Growth{2} \preceq
\Growth{1}$ hold.
\end{definition}

\par The equivalence class of the function $\Growth{}$ is called \textit{the
growth order} and is denoted by the symbol $\GrowthOrder{\Growth{}}$. The
growth order $\GrowthOrder{\Growth{}}$ is called
\begin{enumerate}
\item \textit{polynomial}, if $\GrowthOrder{\Growth{}} = \GrowthOrder{n ^d}$
for some $d > 0$;

\item \textit{intermediate}, if $\GrowthOrder{n ^d} < \GrowthOrder{\Growth{}} <
\GrowthOrder{e ^n}$ for all $d > 0$;

\item \textit{exponential}, if $\GrowthOrder{\Growth{}} = \GrowthOrder{e ^n}$.
\end{enumerate}

\par It is often convenient to encode the growth function of a semigroup in a
generating series:
\begin{definition}
  Let $\Semigroup{}$ be a semigroup generated by a finite set $G$. \emph{The
  growth series} of $\Semigroup{}$ is the formal power series \[
      \Gamma _{\Semigroup{}} \n[X] = \sum \limits _{n \ge 0} \GrowthSemigroup{}
      \n X ^n.
  \]
\end{definition}
The power series $\Delta _{\Semigroup{}} \n[X] = \sum \limits_{n \ge 0}
\WordGrowthSemigroup{} \n X ^n$ can also be introduced; we then have $\Delta
_{\Semigroup{}} \n[X] = (1 - X) \Gamma _{\Semigroup{}} \n[X]$. The series
$\Delta _{\Semigroup{}}$ is called the \textit{word growth series} of the
semigroup $\Semigroup{}$.

\par The growth series of a Mealy automaton is introduced similarly:
\begin{definition}
  Let $\Automaton$ be an arbitrary Mealy automaton. The \emph{growth
    series} of $\Automaton$ is the formal power series \[ \Gamma
  _{\Automaton} \n[X] = \sum \limits _{n \ge 0} \GrowthAutomaton{} \n
  X ^n.
  \]
\end{definition}

\section{Composite growth functions} \label{sect:composite_function}

\par Let us introduce the concept of composite growth function in the following
way. Let $\Growth{} : \Natural \to \Natural$ be an arbitrary function, and let
$k \ge 1$ be a positive integer. Let us define the functions $\Growth{i} :
\Natural \to \Natural$, $i = 0, 1, \ldots, {k - 1}$, by the equalities:
\[
    \Growth{i} \n = \Growth{} \n[k \cdot n + i], \PositiveInteger.
\]
We say that the function $\Growth{}$ is \textit{composite}, if there exists
integer $k \ge 2$ such that at least two functions from the set
\[
    \left\{ \Growth{0}, \Growth{1}, \ldots, \Growth{k - 1} \right\}
\]
can be defined by different expressions.


\par Let us fix the notions. Let $\Automaton$ be the arbitrary Mealy automaton.
Let us denote the semigroup of automatic transformations, defined by
$\Automaton$, by the symbol $\Semigroup{\Automaton}$, and the growth functions
of $\Automaton$ and $\Semigroup{\Automaton}$ by the symbols
$\Growth{\Automaton}$ and $\GrowthSemigroup{\Automaton}$, respectively. If
$\Growth{\Automaton}$ is a composite function for some integer $k$, then let us
denote its ``parts'' by the symbols $\Growth{\Automaton, i}$, $i = 0, 1,
\ldots, k - 1$. Let $\Growth{}$ be an arbitrary function, and let us denote the
$i$-th finite difference of $\Growth{}$ by the symbols $\Growth{} ^{\n[i]}$, $i
\ge 1$, i.e.
\begin{align*}
    \Growth{} ^{\n[1]} \n & = \Growth{} \n - \Growth{} \n[n - 1],\\
    \Growth{} ^{\n[i]} \n & = \Growth{} ^{\n[i - 1]} \n - \Growth{}  ^{\n[i -
    1]} \n[n - 1],
\end{align*}
where $i \ge 2$, $n \ge {i + 1}$.

\begin{figure}[t]
  \centering
  \includegraphics*{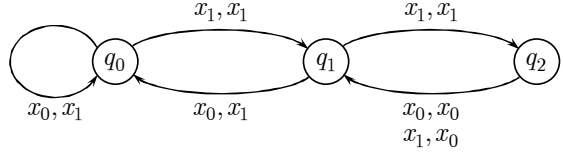}
  \caption{The automaton $\Automaton[1]$}
  \label{fig:automaton_composite exponent}
\end{figure}

\par Let us consider the example of Mealy automaton with the composite growth
function. Let $\Automaton[1]$ be the $3$-state Mealy automaton over the
$2$-symbol alphabet whose Moore diagram is shown on
Figure~\ref{fig:automaton_composite exponent}. Its automatic transformations
have the following unrolled forms:
\begin{align*} \label{eq:unrolled_form_automaton_1}
    f_0 &= \left( f_0, f_1 \right) \left( x_1, x_1 \right),
    & f_1 &= \left( f_0, f_2 \right) \left( x_1, x_1 \right), &
    f_1 &= \left( f_1, f_1 \right) \left( x_0, x_0 \right).
\end{align*}
The following theorem holds:

\begin{theorem} \label{th:semigroup_automaton_1}
\begin{enumerate}
\item The semigroup $\Semigroup{\Automaton[1]}$ has the following presentation:
\begin{equation*} \label{eq:semigroup_automaton_1}
    \Semigroup{\Automaton[5]} = \ATSemigroup
    {\begin{array}{l}
        {f_1 f_2 = f_0 f_2, \, f_1 ^2 f_i = f_1 f_0 f_2,
        \, i = 0, 1; } \\
        {f_0 ^2 f_1 = f_0 ^2, \, f_2 f_0 f_1 = f_2 ^2 f_j = f_2 f_0 ^2, \, j =
        0, 1, 2;} \\
        {f_0 ^4 = f_0 ^3, \, f_0 ^3 f_2 = f_0 ^3, \, f_0 f_2 f_0 ^2 = f_0 f_1
        f_0 f_2;} \\
        {f_1 f_0 f_2 f_0 = f_1 f_0 f_2, \, f_2 f_0 ^3 = f_2 f_0 ^2 }
    \end{array}}
\end{equation*}

\item The growth function $\GrowthAutomaton{1}$ is a composite function for $k
= 2$, and is defined by the following equalities:
\begin{align*}
    \Growth{\Automaton[1], 0} \n & = 23 \cdot 2 ^{n - 2} - 1, &
    \Growth{\Automaton[1], 1} \n & = 32 \cdot 2 ^{n - 2} - 1,
\end{align*}
where $n \ge 2$, $\GrowthAutomaton{1} \n[1] = 3$, $\GrowthAutomaton{1} \n[2] =
8$, $\GrowthAutomaton{1} \n[3] = 14$.
\end{enumerate}
\end{theorem}

\par It follows from Theorem~\ref{th:semigroup_automaton_1} that the growth
function $\GrowthAutomaton{1}$ has the exponential growth order and can be
written as
\begin{equation*}
    \GrowthAutomaton{1} \n =
    \left\{%
    \begin{array}{ll}
        23 \cdot 2 ^{\frac{n - 4}{2}} - 1, & \hbox{if $n$ is even;} \\
        32 \cdot 2 ^{\frac{n - 5}{2}} - 1, & \hbox{if $n$ is odd;} \\
    \end{array}%
    \right.
\end{equation*}
where $n \ge 4$. The normal form of elements of $\Semigroup{\Automaton[1]}$ is
declared in the following proposition:
\begin{proposition} \label{prop:automaton_1_normal_form}
An arbitrary element $\s$ of $\Semigroup{\Automaton[1]}$ has the following
normal form
\begin{equation*}
    \s = \s ' \cdot \left( f_0 f_2 \right) ^{p_1} \left( f_1 f_0 \right) ^{p_2}
    \left( f_0 f_2 \right) ^{p_3} \left( f_1 f_0 \right) ^{p_4} \ldots
    \left( f_0 f_2 \right) ^{p_{2k - 1}} \left( f_1 f_0 \right) ^{p_{2k}} \cdot
    \s '',
\end{equation*}
where $\s' \in \left\{ 1, f_0, f_2 \right\}$, $\s'' \in \left\{ 1, f_0, f_1,
f_1 ^2, f_1 f_0 ^3, f_0 f_2 ^2, f_0 f_2 f_1 f_0 f_2 \right\}$, and $k \,\ge\,
1$, $p_1, p_{2k} \ge 0$, $p_i > 0$, $i = 2, 3, \ldots, {2k - 1}$, $\ell (\s)
\ge 1$.
\end{proposition}

\section{Non-monotonic growth functions} \label{sect:non-monotonic}

\par The conception of a composite function allows to construct
non-monotonic functions easily. For example, let $k = 2$ and $\Growth{}$ be a
function such that $\Growth{1} \n = 1$ and $\Growth{2} \n = 2$. Obviously,
$\Growth{}$ is non-monotonic. Below we consider the $2$-state Mealy automata
over the $4$-symbol alphabet, that have non-monotonic growth functions of
constant, linear and square growth. There are exist automata that have
non-monotonic growth functions of other polynomial growth orders, but their
consideration requires more technical details.

\subsection{The automaton $\Automaton[2]$ of constant growth}

\begin{figure}[t]
  \centering
  \includegraphics*{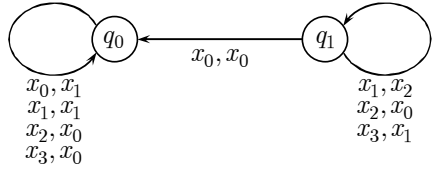}
  \caption{The automaton $\Automaton[2]$}
  \label{fig:automaton_non-monotonic_constant}
\end{figure}

\begin{figure}[b]
  \centering
  \includegraphics*{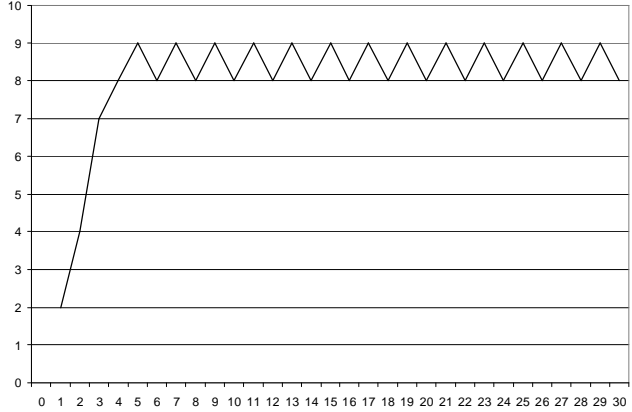}
  \caption{The growth function of $\Automaton[2]$}
  \label{fig:growth_non-monotonic_constant}
\end{figure}

\par Let $\Automaton[2]$ be the Mealy automaton, defined by Moore diagram on
Figure~\ref{fig:automaton_non-monotonic_constant}. Its automatic
transformations have the following unrolled forms:
\begin{align*} \label{eq:unrolled_form_automaton_2}
    f_0 &= \left( f_0, f_0, f_0, f_0 \right) \left( x_1, x_1, x_0, x_0 \right),
    & f_1 &= \left( f_0, f_1, f_1, f_1 \right) \left( x_0, x_2, x_0, x_1
    \right).
\end{align*}
The automaton $\Automaton[2]$ has non-monotonic growth function of the constant
growth order, and the graph of $\GrowthAutomaton{2}$ is shown on
Figure~\ref{fig:growth_non-monotonic_constant}. The following theorem holds:

\begin{theorem} \label{th:semigroup_automaton_2}
\begin{enumerate}
\item The semigroup $\Semigroup{\Automaton[2]}$ has the following presentation:
\begin{equation*} \label{eq:semigroup_automaton_2}
    \Semigroup{\Automaton[2]} = \ATSemigroup
    {\begin{array}{l}
        {f_0 ^2 f_i = f_0 f_1 ^2 f_i = f_0 ^2, \, i = 0, 1; f_1 ^2
        f_0 ^2 = f_0 f_1 f_0 ^2, } \\
        {f_1 f_0 f_1 f_0 ^2 = f_1 ^4 = f_1 ^3 f_0, \, \left( f_1 f_0 \right) ^4
        = \left( f_1 f_0 \right) ^2}
    \end{array}}
\end{equation*}

\item The growth function $\GrowthAutomaton{2}$ is a composite function for $k
= 2$, and is defined by the following equalities:
\begin{align*}
    \Growth{\Automaton[2], 0} \n & = 8, & \Growth{\Automaton[2], 1} \n & = 9,
\end{align*}
where $n \ge 2$, $\GrowthAutomaton{2} \n[1] = 2$, $\GrowthAutomaton{2} \n[2] =
4$, $\GrowthAutomaton{2} \n[3] = 7$.

\end{enumerate}
\end{theorem}

\subsection{The automaton $\Automaton[3]$ of linear growth}

\begin{figure}[t]
  \centering
  \includegraphics*{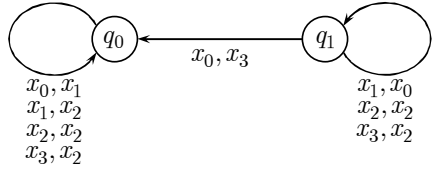}
  \caption{The automaton $\Automaton[3]$}
  \label{fig:automaton_non-monotonic_linear}
\end{figure}

\begin{figure}[b]
  \centering
  \includegraphics*{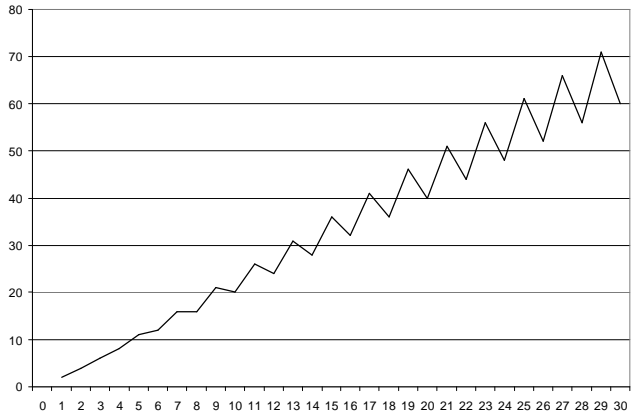}
  \caption{The growth function of $\Automaton[3]$}
  \label{fig:growth_non-monotonic_linear}
\end{figure}

\par Let's consider the automaton $\Automaton[3]$, whose Moore diagram is shown on
Figure~\ref{fig:automaton_non-monotonic_linear}. Its automatic transformations
have the following unrolled forms:
\begin{align*} \label{eq:unrolled_form_automaton_3}
    f_0 &= \left( f_0, f_0, f_0, f_0 \right) \left( x_1, x_2, x_2, x_2 \right),
    & f_1 &= \left( f_0, f_1, f_1, f_1 \right) \left( x_3, x_0, x_2, x_2
    \right).
\end{align*}
The automaton $\Automaton[3]$ have the non-monotonic linear growth function,
and the graph of $\GrowthAutomaton{3}$ is shown on
Figure~\ref{fig:growth_non-monotonic_linear}. The following theorem holds:

\begin{theorem} \label{th:semigroup_automaton_3}
\begin{enumerate}
\item The semigroup $\Semigroup{\Automaton[3]}$ has the following presentation:
\begin{equation*}
    \Semigroup{\Automaton[3]} = \ATSemigroup
    {\begin{array}{l}
        {f_0 ^2 f_i = f_1 f_0 ^2 = f_0 ^2, \, i = 0, 1; f_0 f_1 f_0 = f_0, } \\
        {f_0 f_1 ^2 f_0 = f_0 ^2, \, f_1 ^3 f_0 f_1 = f_1 f_0 f_1 ^3}
    \end{array}}
\end{equation*}

\item The growth function $\GrowthAutomaton{3}$ is a composite function for $k
= 2$, and is defined by the following equalities:
\begin{align*}
    \Growth{\Automaton[3], 0} \n & = 4n, & \Growth{\Automaton[3], 1} \n & = 5n
    + 1,
\end{align*}
where $n \ge 1$ and $\GrowthAutomaton{3} \n[1] = 2$.

\end{enumerate}
\end{theorem}

\subsection{The automaton $\Automaton[4]$ of square growth}

\begin{figure}[t]
  \centering
  \includegraphics*{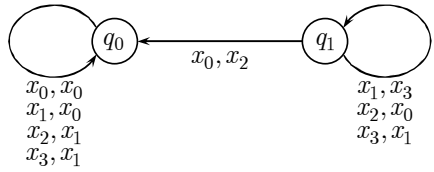}
  \caption{The automaton $\Automaton[4]$}
  \label{fig:automaton_non-monotonic_square}
\end{figure}

\begin{figure}[b]
  \centering
  \includegraphics*{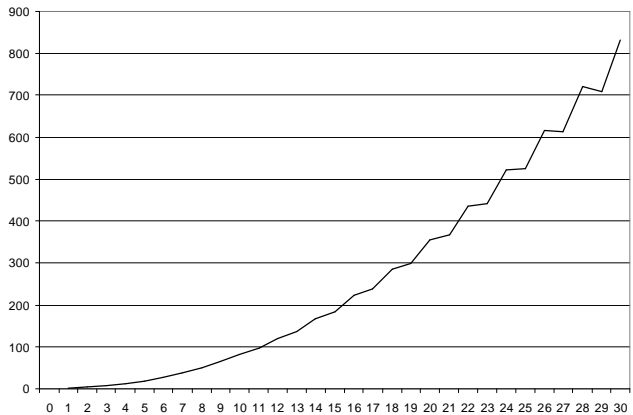}
  \caption{The growth function of $\Automaton[4]$}
  \label{fig:growth_non-monotonic_square}
\end{figure}

\par Let $\Automaton[4]$ be the automaton such that its Moore diagram is shown on
Figure~\ref{fig:automaton_non-monotonic_square}. Its automatic transformations
have the following unrolled forms:
\begin{align*} \label{eq:unrolled_form_automaton_4}
    f_0 &= \left( f_0, f_0, f_0, f_0 \right) \left( x_0, x_0, x_1, x_1 \right),
    & f_1 &= \left( f_0, f_1, f_1, f_1 \right) \left( x_2, x_3, x_0, x_1
    \right).
\end{align*}

\par The automaton $\Automaton[4]$ has non-monotonic growth function of square growth, and the
graph of $\GrowthAutomaton{4}$ is shown on
Figure~\ref{fig:growth_non-monotonic_square}. The properties of $\Automaton[4]$
are formulated in the following theorem:

\begin{theorem} \label{th:semigroup_automaton_4}
\begin{enumerate}
\item The semigroup $\Semigroup{\Automaton[4]}$ is infinitely presented, and
has the following presentation:
\begin{equation} \label{eq:semigroup_automaton_4}
    \Semigroup{\Automaton[4]} = \ATSemigroup
    {\begin{array}{l}
        {f_0 ^2 f_i = f_0 ^2, \, f_1 f_0 f_1 ^2 f_i = f_1 f_0 f_i, \, i = 0, 1,
        } \\
        {f_0 f_1 ^{2p + 1} f_0 = f_0 f_1 f_0, \, p \ge 1}
    \end{array}}
\end{equation}

\item The growth function $\GrowthAutomaton{4}$ is a composite non-monotonic
function, that is defined by the following equalities:
\begin{align*}
    \Growth{\Automaton[4], 0} \n & = 4n ^2 - 5n + 6, n \ge 2, & \Growth{\Automaton[4],
    1} \n & = \frac{7}{2}n ^2 + \frac{3}{2}n + 2, n \ge 0
\end{align*}
and $\GrowthAutomaton{4} \n[2] = 4$. The function $\GrowthAutomaton{4}$ has the
square growth order.

\end{enumerate}
\end{theorem}

\par From the defining relations~\eqref{eq:semigroup_automaton_4} the
proposition follows
\begin{proposition} \label{prop:automaton_4_normal_form}
An arbitrary element $\s$ of $\Semigroup{\Automaton[4]}$ admits a unique
minimal-length representation as a word of one of the following forms
\begin{equation*}
    f_0 f_1 ^{2p_1} \left( f_0 f_1 \right) ^{p_2} \cdot \s',
\end{equation*}
where $p_1 \ge 1$, $p_2 \ge 0$, $\s' \in \left\{ 1, f_1, f_0, f_0 ^2 \right\}$,
except the combination $p_1 = 1$, $p_2 = 0$, $\s' = 1$, or
\begin{equation*}
    f_1 ^{p_1} \left( f_0 f_1 \right) ^{p_2} \cdot \s',
\end{equation*}
where $p_1 \ge 0$, $p_2 \ge 1$, $\s' \in \left\{ 1, f_1, f_0, f_0 ^2 \right\}$,
or $p_2 = 0$, $p_1 \ge 0$, $\s' \in \left\{ f_1, f_0, f_0 ^2 \right\}$.
\end{proposition}

\section{Growth functions with doubled finite differences} \label{sect:doubled_differnces}

\par In this section we consider composite growth functions such that one of
its finite differences consists of doubled values. We consider the sequence
$\SequenceIndex{B}{3}$ of Mealy automata of polynomial growth, and two Mealy
automata of the intermediate and the exponential growth orders.

\subsection{The automata $\SequenceIndex{B}{3}$ of polynomial growth}

\begin{figure}[t]
  \centering
  \includegraphics*{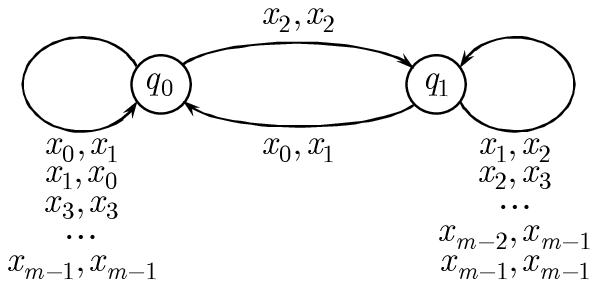}
  \caption{The automaton $\AExample{m}$}
  \label{fig:automaton_polynomial}
\end{figure}

\par Let $\AExample{m}$, $m \ge 3$, be the $2$-state Mealy automaton over the
$m$-symbol alphabet (Figure~\ref{fig:automaton_polynomial}), and the unrolled
forms of the automatic transformations $f_0 = f_{q_0, \AExample{m}}$ and $f_1 =
f_{q_1, \AExample{m}}$ are defined in the following way:
\begin{align*}
    f_0 & = \left( f_0, f_0, f_1, f_0, \ldots, f_0, f_0 \right)
    \left( x_1, x_0, x_2, x_3, \ldots, x_{m - 2}, x_{m - 1} \right),\\
    f_1 & = \left( f_0, f_1, f_1, f_1, \ldots, f_1, f_1 \right)
    \left( x_1, x_2, x_3, x_4, \ldots, x_{m - 1}, x_{m - 1} \right).
\end{align*}

\begin{theorem}  \label{th:semigroups_polynomial}
\begin{enumerate}
\item For any $m \ge 3$ the semigroup $\Semigroup{m}$ have the following
presentation:
\begin{align*}
    \Semigroup{\AExample{3}} & = \ATSemigroup{
        f_1 ^3 = f_0 f_1 ^2, f_1 f_0 f_1 = f_0 ^2 f_1},\\
    \Semigroup{\AExample{4}} & = \ATSemigroup{
        f_1 ^4 = f_1 f_0 f_1 ^2, \, f_1 f_0 ^{p_1} f_1 f_0 f_1 = f_1 f_0 ^{p_1
        + 2} f_1, p_1 \ge 0 },\\
    \Semigroup{\AExample{m}} & = \ATSemigroup{
        {\begin{array}{l}
            {\prod \limits _{i = 1} ^{m - 4}{\left( { f_1 f_0 ^{p_i}} \right)} f_1
            ^4 = \prod \limits _{i = 1} ^{m - 4}{\left( { f_1 f_0 ^{p_i}} \right)}
            f_1 f_0 f_1 ^2,}\\
            {\prod \limits _{i = 1} ^{m - 3}{\left( { f_1 f_0 ^{p_i}} \right)} f_1
            f_0 f_1 = \prod \limits _{i = 1} ^{m - 3}{\left( { f_1 f_0 ^{p_i}} \right)}
            f_0 ^2 f_1,}\\
            {p_i \ge 0, i = 1, 2, \ldots, {m - 3}}
        \end{array}}}.
\end{align*}
All semigroups $\Semigroup{\AExample{m}}$ for $m \ge 4$ are infinitely
presented.

\item For $m \ge 3$ the growth function $\Growth{\AExample{m}}$ is defined by
the following equalities:
\begin{multline} \label{eq:growth_function_aexample_m}
    \Growth{\AExample{m}} \n = \sum \limits _{i = 0} ^{m - 2} {\Binomial{i}{n}}
    + \sum \limits _{i = 0} ^{\Divider{n - m + 1}} {\Binomial{m - 2}{n - 2i -
    1}} =\\
    = \sum \limits _{i = 0} ^{m - 2} {\Binomial{i}{n}} + \sum \limits _{i \ge
    0} {\Binomial{m - 2}{n - 2i - 1}},
\end{multline}
for all $n \ge 1$.

\end{enumerate}
\end{theorem}

\par Here $\left[ r \right]$ denotes the integer part of the real number $r$,
and we assume that $\Binomial{k}{n} = 0$ if $k \ge n$ or $n < 0$.
\par The following proposition holds in the semigroup $\Semigroup{\AExample{m}}$:
\begin{proposition}
The normal form of the element $\s$ of $\Semigroup{\AExample{m}}$ is one of the
following words
\begin{equation*}
    f_0 ^{p_1} f_1 f_0 ^{p_2} f_1 \ldots f_0 ^{p_{k - 1}} f_1 f_0 ^{p_k}
\end{equation*}
where $1 \le k \le {m - 1}$, $p_i \ge 0$, $i = 1, 2, \ldots, {k}$, $\ell(\s)
\ge 1$, and
\begin{equation*}
    f_0 ^{p_1} f_1 f_0 ^{p_2} f_1 \ldots f_0 ^{p_{m - 2}} f_1 f_0 ^{2p_{m - 1}}
    f_1 f_0 ^{p_m}
\end{equation*}
where $p_i \ge 0$, $i = 1, 2, \ldots, {m}$.
\end{proposition}

\par The corollary follows from
Theorem~\ref{th:semigroups_polynomial}:
\begin{corollary}
\begin{enumerate}
\item For all $m \ge 3$ the function $\Growth{\AExample{m}}$ have the growth
order $\GrowthOrder{n^{m - 1}}$.

\item The $\left( m - 2 \right)$-th finite differences of
$\Growth{\AExample{m}}$ is defined by the following equality
\[
    \Growth{\AExample{m}} ^{\n[m - 2]} \n = \Divider{n - m + 1} + 2,
\]
where $n \ge m - 1$.
\end{enumerate}
\end{corollary}

\par It follows from~\eqref{eq:growth_function_aexample_m} that for any $m \ge
4$ the equalities hold
\begin{multline*}
    \Growth{\AExample{m}} ^{\n[1]} \n = \Growth{\AExample{m}} \n -
    \Growth{\AExample{m}} \n[n - 1] = \\
    = \sum \limits _{i = 0} ^{m - 3} {\Binomial{i}{n - 1}} + \sum \limits _{i
    \ge 0} {\Binomial{m - 3}{n - 2i - 2}} = \Growth{\AExample{m - 1}} \n[n -
    1],
\end{multline*}
where $n \ge 2$. The growth function $\Growth{\AExample{3}}$ is defined by the
equalities
\begin{equation*}
    \Growth{\AExample{3}} \n =
    \left\{%
    \begin{array}{ll}
        \frac{1}{4}n ^2 + n + 1, & \hbox{if $n$ is even;} \\
        \frac{1}{4}n ^2 + n + \frac{3}{4}, & \hbox{if $n$ is odd.} \\
    \end{array}%
    \right.
\end{equation*}
Hence, joining two last equalities, one has
\[
    \Growth{\AExample{m}} ^{\n[m - 2]} \n = \Growth{\AExample{3}} ^{\n[1]} \n[n
    - \left( m - 3 \right)] = \Divider{n - \left( m - 3 \right)} + 1 =
    \Divider{n - m + 1} + 2,
\]
for any $m \ge 3$, $n \ge {m - 1}$. Therefore $\left( m - 2 \right)$-th finite
difference of $\Growth{\AExample{m}}$ consists of doubled values, i.e. for any
even integer $n$, $n \ge {0}$, the equality holds
\[
    \Growth{\AExample{m}} ^{\n[m - 2]} \n[n + m] =
    \Growth{\AExample{m}} ^{\n[m - 2]} \n[n + m - 1].
\]

\subsection{The automaton $\Automaton[5]$ of intermediate growth}

\begin{figure}[t]
  \centering
  \includegraphics*{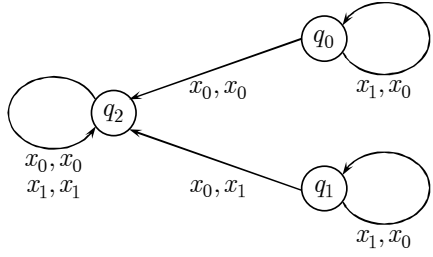}
  \caption{The automaton $\Automaton[5]$}
  \label{fig:automaton_intermediate}
\end{figure}

\par Let $\Automaton[5]$ be the $3$-state Mealy automaton over the
$2$-symbol alphabet such that its Moore diagram is shown on
Figure~\ref{fig:automaton_intermediate}. The following theorem holds:

\begin{theorem} \label{th:semigroup_intermediate}
\begin{enumerate}
\item The semigroup $\Semigroup{\Automaton[5]}$ is a infinitely presented
monoid, and has the following presentation:
\begin{equation} \label{eq:semigroup_automaton_5}
    \Semigroup{\Automaton[5]} = \ATMonoid{
        {\begin{aligned}
            & f_0 f_1 ^{2^k - 1} \cdot f_1 ^ {p 2 ^{k + 1}} f_0 \prod
            \limits_{i = k}^{1} \left( f_1 ^ {2 ^{i} - 1} f_0 \right) =\\
            & = f_1 ^ {p 2 ^{k + 1}} f_0 \prod \limits_{i = k}^{1} \left( f_1 ^
            {2 ^{i} - 1} f_0 \right), k \ge 0, p = 0, 1.
        \end{aligned}}}
\end{equation}

\item The growth series $\Gamma _{\Automaton[5]} \n[X]$ of $\Automaton[5]$ and
the growth series $\Gamma _{\Semigroup{\Automaton[5]}} \n[X]$ of
$\Semigroup{\Automaton[5]}$ coincide and are defined by the equality
\begin{multline*}
    \Gamma _{\Semigroup{\Automaton[5]}} \n[X] = \Gamma _{\Automaton[5]} \n[X] =
    \frac{1}{(1 - X) ^2} \left( 1 + \frac{X}{1 - X} \left( 1 +
    \frac{X^2}{1 - X ^2} \cdot \right. \right. \\
    \left. \left. \left. \left( 1 + \frac{X^4}{1 - X ^4} \left( 1 +
    \frac{X^8}{1 - X ^8} \left(1 + \ldots \right) \right) \right) \right)
    \right) \right)
\end{multline*}

\end{enumerate}
\end{theorem}

\par The properties of the growth function $\GrowthAutomaton{5}$ are formulated
in the following corollary:
\begin{corollary}
\begin{enumerate}
\item The growth function $\GrowthAutomaton{5}$ has the intermediate growth
order $\GrowthOrder{ n ^{\frac{\log n}{2 \log 2}}}$.

\item Let us define $\GrowthAutomaton{5} ^{\n[2]} \n[0] = \GrowthAutomaton{5}
^{\n[2]} \n[1] = \GrowthAutomaton{5} ^{\n[2]} \n[2] = 1$. The second finite
difference of $\GrowthAutomaton{5}$ is defined by the following equality
\begin{equation} \label{eq:finite_differences_intermediate}
    \GrowthAutomaton{5} ^{\n[2]} \n = \sum \limits _{i = 0} ^{\Divider{n - 1}}
    \GrowthAutomaton{5} ^{\n[2]} \n[i], \, n \ge 3.
\end{equation}
\end{enumerate}
\end{corollary}

\par The system of defining relations \eqref{eq:semigroup_automaton_5} implies
the following normal form:
\begin{proposition}
Each semigroup element $\s$ of $\Semigroup{\Automaton[5]}$ can be written in
the following normal form
\begin{equation*}
    f_1 ^{p_0} f_0 f_1 ^{2^{k - 1} p_1 + ({2^{k - 1} - 1})} f_0 \ldots
    f_1 ^{2^{i} p_{k - i} + (2^{i} - 1)} f_0 \ldots
    f_1 ^{4 p_{k - 2} + 3} f_0 f_1 ^{2p_{k - 1} + 1} f_0 f_1 ^{p_k}
\end{equation*}
where $k \ge 0$, $p_i \ge 0$, $i = 0, 1, \ldots, {k}$.
\end{proposition}

\par The growth series for the second finite difference $\Delta ^{\n[2]} \Gamma
_{\Automaton[5]} \n[X]$ can be easily constructed by using the expression for
$\Gamma _{\Automaton[5]} \n[X]$:
\begin{align*}
    \Delta ^{\n[2]} & \Gamma _{\Automaton[5]} \n[X] = \sum \limits _{n \ge 3}
    {\GrowthAutomaton{5} ^{\n[2]} \n X^n} + \GrowthAutomaton{5} ^{\n[2]} \n[0]
    + \GrowthAutomaton{5} ^{\n[2]} \n[1] X + \GrowthAutomaton{5} ^{\n[2]} \n[2]
    X^2 =\\
    & {\begin{aligned}
    = (1 - X)^2 & \Gamma _{\Automaton[5]} \n[X] - (1 - X)
    \GrowthAutomaton{5} \n[0] - X (\GrowthAutomaton{5} \n[1] -
    \GrowthAutomaton{5} \n[0]) - \\
    & - X ^2 (\GrowthAutomaton{5} \n[2] - 2 \GrowthAutomaton{5} \n[1]
    + \GrowthAutomaton{5} \n[0]) + 1 + X + X^2 =
    \end{aligned}}\\
    & = 1 + \frac{X}{1 - X} \left( 1 + \frac{X^ 2}{1 - X ^2} \left( 1 +
    \frac{X^4}{1 - X ^4} \left( 1 + \frac{X^8}{1 - X ^8} \left(1 + \ldots
    \right) \right) \right) \right).
\end{align*}
The right-hand series of the last equality are the formal series for the
numbers of partitions of $n$, $n \ge 1$, into ``sequential'' powers of $2$,
that is $\GrowthAutomaton{5} ^{\n[2]} \n$ equals the cardinality of the set
\[
    \left\{
        \begin{array}{*{20}c}
           {p_0, p_1, \ldots, p_k} & \vline & {k \ge 0, \sum \limits_{i = 0}
           ^{k} p_i 2^i = n, p_i \ge 1, i = 0, 1, \ldots, k} \\
        \end{array}
    \right\}.
\]
The equality~\eqref{eq:finite_differences_intermediate} is well-known for these
partition numbers \cite{Andrews1976}. Therefore, the second finite difference
of $\GrowthAutomaton{5}$ consists of doubled values, i.e. for any even integer
$n$, $n \ge 2$, the equality holds
\[
    \GrowthAutomaton{5} ^{\n[2]} \n = \GrowthAutomaton{5} ^{\n[2]}
    \n[n - 1].
\]

\subsection{The automaton $\Automaton[6]$ of exponential growth}

\begin{figure}[t]
  \centering
  \includegraphics*{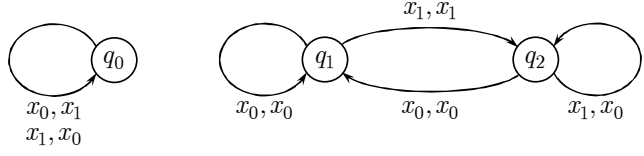}
  \caption{The automaton $\Automaton[6]$}
  \label{fig:automaton_Fibonacci}
\end{figure}

\par Let $\Automaton[6]$ be the $3$-state Mealy automaton over the
$2$-symbol alphabet such that its automatic transformations have the following
unrolled forms:
\begin{align*} \label{eq:unrolled_form_automaton_6}
    f_0 &= \left( f_0, f_0 \right) \left( x_1, x_0 \right), &
    f_1 &= \left( f_1, f_2 \right) \left( x_0, x_1 \right), &
    f_2 &= \left( f_1, f_2 \right) \left( x_0, x_0 \right).
\end{align*}
The Moore diagram of $\Automaton[6]$ is shown on
Figure~\ref{fig:automaton_Fibonacci}. The following theorem holds:

\begin{theorem} \label{th:semigroup_automaton_6}
\begin{enumerate}
\item The semigroup $\Semigroup{\Automaton[6]}$ has the following presentation:
\begin{equation*}
    \Semigroup{\Automaton[6]} = \ATSemigroupThree
    {\begin{array}{l}
        {f_0 ^2 = 1, \, f_2 f_1 = f_1 f_2 = f_2 ^2 = f_2} \\
        {f_1 ^2 = f_1, \, f_2 f_0 f_1 f_0 f_2 = f_1 f_0 f_1 f_0 f_2}
    \end{array}}
\end{equation*}

\item The growth series $\Gamma _{\Automaton[6]} \n[X]$ of $\Automaton[6]$
admits the description
\begin{equation*}
    \Gamma _{\Automaton[6]} \n[X] = \frac{1}{(1 - X) ^2}
    \left( 2 X - 1 + \frac{1 + X + X ^3}{1 - X ^2 - X ^4} \right).
\end{equation*}

\item The growth series $\Gamma _{\Semigroup{\Automaton[6]}} \n[X]$ of
$\Semigroup{\Automaton[6]}$ is defined in the following way
\begin{equation*}
    \Gamma _{\Semigroup{\Automaton[6]}} \n[X] = \frac{1}{(1 -
    X) ^2} \left( X + \frac{1 + X + X ^3}{1 - X ^2 - X ^4} \right).
\end{equation*}

\end{enumerate}
\end{theorem}

\par Let us define the Fibonacci numbers by the symbols $\Phi _n$, where $\Phi _n
= \Phi_ {n - 1} + \Phi_ {n - 2}$, $n \ge 2$, and $\Phi_ {0} = \Phi _{1} = 1$.
It follows from Theorem~\ref{th:semigroup_automaton_6} that the growth function
$\GrowthAutomaton{6}$ can be written in close form, and the following corollary
holds:
\begin{corollary} \label{th:growth_functions_automaton_6}
The growth function $\GrowthAutomaton{6}$ is defined by the following
equalities:
\begin{equation} \label{eq:growth_function_automaton_10}
    \Growth{\Automaton[6]} \n =
    \left\{%
    \begin{array}{ll}
        \Phi_{\Divider{n} + 6} + \Phi_{\Divider{n} + 4} - 2n - 18, & \hbox{if
        $n$ is even;} \\
        \Phi_{\Divider{n} + 6} + 2 \Phi_{\Divider{n} + 4} - 2n - 18, &
        \hbox{if $n$ is odd.} \\
    \end{array}%
    \right.
\end{equation}
The growth function $\GrowthAutomaton{6}$ has the exponential growth order.
\end{corollary}

\par Let $n$ be any positive integer, and represent $n = 2k$, when $n$ is even,
and $n = 2k + 1$, when $n$ is odd. It follows from
\eqref{eq:growth_function_automaton_10}, that for any $k \ge 0$ the following
equalities hold
\begin{align*}
    \Growth{\Automaton[6]} ^{\n[1]} \n[2k + 1] & = \Phi_{k + 4} - 2, &
    \Growth{\Automaton[6]} ^{\n[1]} \n[2k + 2] & = 2 \Phi_{k + 3} - 2,
\end{align*}
and, using the previous equalities, we have
\begin{align*}
    \Growth{\Automaton[6]} ^{\n[2]} \n[2k + 1] & = \Phi_{k + 1}, &
    \Growth{\Automaton[6]} ^{\n[2]} \n[2k + 2] & = \Phi_{k + 1}.
\end{align*}
Hence, the second finite difference $\Growth{\Automaton[6]}$ consists of
doubled values, and for all even integer $n$ the equality holds
\[
    \Growth{\Automaton[6]} ^{\n[2]} \n[n] = \Growth{\Automaton[6]} ^{\n[2]}
    \n[n - 1].
\]

\section{Final remarks} \label{sect:final_remarks}

\par There are some questions, that concern the composite non-monotonic growth
functions of Mealy automata.

\begin{enumerate}
\item Does there exist the Mealy automaton such that its composite growth
function includes ``parts'' of different growth orders?

\item Does there exist the Mealy automaton which have the non-monotonic growth
function of the intermediate or the exponential growth order?

\item Does there exist the Mealy automaton such that its growth function is
non-monotonic, but isn't a composite function (in the sense of
Section~\ref{sect:composite_function})?

\end{enumerate}

\end{document}